\documentclass[11pt]{amsart}
\usepackage{amscd,amssymb}

\usepackage[cmtip,matrix,arrow]{xy}

\numberwithin{equation}{section}

\newtheorem {Theorem} 			{Theorem}

\newtheorem {Proposition}[equation]	{Proposition}
\newtheorem {Conjecture}[equation]	{Conjecture}

\theoremstyle{definition}

\newtheorem {Remark}[equation]		{Remark}

\newcommand{\pr} {\smallskip\noindent{\bf Proof\,\,}}

\newcommand     {\comment}[1]   {}
\newcommand     {\mute}[2] {}
\newcommand     {\printname}[1] {}

\newcommand{\labell}[1] {\label{#1}\printname{#1}}

\def	\to	{\longrightarrow}

\def	\C	{{\mathbb C}}
\def	\R	{{\mathbb R}}

\def	\Z	{{\mathbb Z}}
\def	\P	{\mathbb P}
\def\cA {{\mathcal A}}
\def\cG {{\mathcal G}}

\def	\cP	{{\mathcal P}}
\def	\cF	{{\mathcal F}}
\def	\cL	{{\mathcal L}}

\def	\fg	{{\mathfrak g}}

\def	\pr	{\operatorname{pr}}

\def\l{\lambda}
\def\G{\Gamma}
\def\a{\alpha}
\def\cB{\mathcal B}
\def\cS{\mathcal S}

\def\tcg{\tilde{\cG}}

\def\cH{\mathcal H}
\def\cC{\mathcal C}

\begin{document}


\title[Measures on Banach manifolds]
{Measures on Banach Manifolds and Supersymmetric Quantum Field Theory}
\author{Jonathan Weitsman}
\thanks{Supported in part by NSF grant DMS 04/05670}
\thanks\today
\address{Department of Mathematics, University of California, Santa Cruz,
CA 95064}
\email{weitsman@math.ucsc.edu}
\begin{abstract}
We show how to construct measures on Banach manifolds associated
to supersymmetric quantum field theories.   These measures are mathematically
well-defined objects inspired by
the formal path integrals appearing in the physics literature
on quantum field theory.  We give three concrete
examples of our construction.  The first example is a family $\mu_P^{s,t}$
of measures
on a space of functions on the two-torus, parametrized by a 
polynomial $P$ (the Wess-Zumino-Landau-Ginzburg model).
The second is a family $\mu_\cG^{s,t}$ of measures on a space $\cG$ of 
maps from $\P^1$ to a Lie group
(the Wess-Zumino-Novikov-Witten model).  Finally we study
a family  $\mu_{M,G}^{s,t}$ of measures on the product of a space of connections on the
trivial principal bundle with structure group $G$ on a three-dimensional
manifold $M$ with a space of $\fg$-valued three-forms on $M.$

We show that these measures are positive, and that the measures $\mu_\cG^{s,t}$
are Borel probability measures.
As an application we show that formulas arising from expectations
in the measures $\mu_\cG^{s,1}$ reproduce formulas discovered by Frenkel and Zhu in
the theory
of vertex operator algebras.  We conjecture that a similar computation for
the measures $\mu_{M,SU(2)}^{s,t},$ where $M$ is a homology three-sphere,
will yield the Casson invariant of $M.$
\end{abstract}
\maketitle

\section{Introduction}

During the past two decades techniques and ideas arising from
quantum field theory have played an increasing role in geometry
and topology.  The work of the school of mathematical
physicists led by E. Witten on applications of supersymmetric
quantum field theory and later of string theory to problems
in mathematics has in many cases led to results that can be stated
as mathematically precise conjectures, and often rigorously proved as
theorems.  However, the methods used by the physicists, among which
formal path-integral methods of integration over
infinite-dimensional spaces are prominent, have not in the main
been accessible to mathematical interpretation; for actual computation
one proceeds by inspired analogy to finite-dimensional cases
and by physical intuition.  The purpose of this paper is
to produce measures on infinite-dimensional spaces associated
to supersymmetric quantum field theories which maintain the
formal analogies to finite dimensions used by the physicists,
but which are mathematically well-defined.

\subsection{Finite dimensional analogs of supersymmetric quantum field theory}
\labell{finsec}

In this paper we give a construction of measures on
Banach manifolds associated to supersymmetric quantum field theories.
In order to approach this problem, we begin with finite-dimensional
analogs of these theories and the associated measures.
In the simplest case, given a map of
finite-dimensional vector spaces, we give three constructions that
coincide:  That of the degree of the map, of the pull-back by this
map of a Gaussian differential form of top degree, and of a measure
obtained by pushforwards by local inverses.  The third construction
is the one which we generalize to Banach manifolds.  The physics literature
makes use of formal analogy to finite dimensions to work with
the first and second constructions only.

\subsubsection{Measures arising from maps of vector spaces}
Let $V,W$ be real, $n$-dimensional Hermitian vector spaces.
Let $\eta_W$ denote the invariant volume form on $W$, and let
$\xi \in H^n(W)$ denote the differential form 
$\xi = \frac{1}{{(2\pi)}^{n/2}} \exp (-|w|^2/2) \eta_W.$  Given a proper
smooth 
map $f : V \to W,$ the pullback $f^*\xi$ gives a signed measure
on $V.$  In terms of the natural volume form $\eta_V$ this measure
is given by 
\begin{equation}\labell{pb}
f^*\xi = \frac{1}{{(2\pi)}^{n/2}} 
\exp (-|f|^2/2) {\rm det}(\nabla f) \eta_V.
\end{equation}

If $f$ grows rapidly at infinity the integral of $f^*\xi$ converges and gives
the Leray-Schauder degree of the map $f,$ so that if $0$
is a regular value of $f$, we have

\begin{equation}\labell{ls}
\int_V f^*\xi = \sum_{x \in f^{-1}(0)} \epsilon_x
\end{equation}

\noindent where $\epsilon_x={\rm det}(\nabla f|_x)/|{\rm det}(\nabla f|_x)|.$
Equation (\ref{ls}) can be seen
to follow formally from the
expression (\ref{pb}) by replacing $f$ by $t f$ and then letting $t \to 
\infty.$  The determinant appearing in (\ref{pb}) can be viewed
as arising from a Berezin integration, so that the integral
is a finite-dimensional analog of supersymmetric quantum field
theory.  The fact that supersymmetric field theories are analogs
of integrals of the type (\ref{pb}) and are related to the
degree of the map corresponding to $f$ is due to 
Witten in \cite{tr-1f,witmor}.  

Formulas similar to (\ref{pb}) and (\ref{ls}) hold
where $V,W$ are complex Hermitian vector spaces and $f: V \to W$
is a proper smooth map.
If $W$ is a  complex vector space of complex dimension $n,$
the Euclidean volume form
on $W$ can be written as $\eta_W = \omega^n/n!,$
where $\omega$ is a two-form on $W$ given in terms of a set 
$(w_1,\dots,w_n)$ of holomorphic coordinates on $W$
by $\omega = \sum_i dw_i \wedge d{\bar{w_i}}.$
It follows that if
$f:V\to W$ is holomorphic and $x \in V$ is a regular point
of $f$ then ${\rm det}(\nabla f|_x) /|{\rm det}(\nabla f|_x)| = 1.$ 

An expression for the measure given by $f^*\xi$ can also be given
as follows.  Let $x \in V$ be a regular point of $f,$ and let $U$
be a neighborhood of $x$ where $f$ is invertible; denote this local
inverse by $g$.  Let $\mu_W$ be the Gaussian measure on $W$ corresponding
to $\xi,$ and consider its restriction to $f(U)$.
We define a measure on $U$ by 

\begin{equation}\labell{measdef}
\mu_U := g_*( \mu_W|_{f(U)}).
\end{equation}

This gives a measure $f^*\mu_W$ on $V$ by taking the union over all neighborhoods
$U$ and by assigning measure zero to the singular points of $f.$
The measure $f^*\mu_W$ can be obtained from the signed measure associated
to $f^*\xi$ by taking absolute values of the determinant of $\nabla f;$
in other words, for any $\phi \in C_c(V)$,

\begin{equation} \labell{phasofmeas}
\frac{1}{{(2\pi)}^{n/2}} \int \phi
\exp (-|f|^2/2) {\rm det}(\nabla f) \eta_V = 
f^*\mu_W (\phi {\rm ~Phase~ det} (\nabla f))
\end{equation}

\noindent where  ${\rm Phase ~det}(\nabla f)(x)$ is given
by ${\rm Phase ~det}(\nabla f)(x): =  {\rm det}(\nabla f|_x)/|{\rm det}(\nabla f_x)|$ if $x$ is
a regular point of $f$ and is extended arbitrarily
to the singular set of $f$ where the measure $f^*\mu_W$ is
zero.

Two features of the finite-dimensional case will recur in
the generalization of this construction to Banach manifolds.
These are the subtlety involved in the definition of the
phase of the determinant of $\nabla f,$ and the issue of
the finiteness of the measure $f^* \mu_W.$
In the case of Banach manifolds,
the measure analogous to $f^*\mu_W$ will have a straightforward
construction, but the phase of the determinant 
may not be well-defined.  This phenomenon is known in
the physics literature as the appearance of an anomaly; the
geometry underlying anomalies was studied by Atiyah-Singer \cite{as}
and Quillen \cite{q}. Similarly, in finite
dimensions, the finiteness of the measure $f^*\mu_W$ depends on the
growth properties of $f$ at infinity, and must be studied
by hard analysis.  The same will be true in the case of Banach
manifolds.

\subsubsection{Measures arising from sections of vector bundles}
Finite-dimensional analogs of supersymmetric quantum field theories
also arise in a slightly different context.  Let $M^n$ be a compact Riemannian
manifold, let $E\to M$ be an $n$-dimensional Hermitian vector bundle over $M$, 
and let $\nabla$ be a connection on $E.$
Let $s: M \to E$ be a  section of $E.$  In this case Mathai and Quillen
\cite{mq} give an expression for the Euler number $\chi(E)$ of $E$ as 

\begin{equation}\labell{mqf}
\chi(E) = \frac{1}{{(2\pi)}^{n/2}} \int_M \exp(-|s|^2/2)
{\rm det}(\nabla s) \sum_{k=1}^n C_k;
\end{equation}

\noindent here $C_k$ is a sequence of terms involving
$\nabla s$ and 
and the curvature of $E.$  In terms of local 
coordinates the Mathai-Quillen formula is given by the integral
of a sum of such terms multiplied by a phase of a determinant
with respect to a measure of the type defined in the previous section.

\subsection{Outline of the paper}

The physics literature on topological quantum field theory is concerned
with formal infinite-dimensional 
integrals corresponding morally to the differential form (\ref{pb}), where 
the map $f:V \to W$ is replaced by a map
$\cF$ of infinite dimensional spaces. It then seeks to apply
ideas from differential geometry (such as localization), as well
as novel techniques, to gain insight
about the solutions of the equation $\cF = 0,$ and more generally about
the formal properties of the underlying measure.  While some of the results
of this approach can be stated mathematically and sometimes verified
by various mathematical methods, there have been so far no satisfactory
constructions of measures on infinite dimensional spaces which would
approach the physicists' formal constructions.

In this paper we will describe a method of constructing measures
on Banach manifolds
analogous to the pullback measures $f^*\mu_W$ in a number
of examples motivated by supersymmetric quantum field theories.
Starting with a supersymmetric model based on a map $\cF$ of
infinite-dimensional spaces, we
will find a family of Banach manifolds $X_s, s > s_0,$ a Banach space $Y,$
and a family of maps
$\cF_s: X_s \to Y,$ where the set $\cF_s = 0$ is independent of $s$ 
and matches formally the set $\cF_0=\cF=0.$  The Banach
space $Y$ in our examples will be equipped with a white noise
measure---the infinite-dimensional analog of the Gaussian measure
$\mu_W$ above.  We will construct an analog of the pullback
measure $f^* \mu_W$ and of the phase of the determinant.  We do
not know of a general theorem that would give convergence of the
integral of the phase of the determinant; instead we study
some examples paralleling supersymmetric quantum
field theories in the literature.  In these examples we can give
specific conjectures about the finiteness of the measure and the
value of the integral of the phase of the determinant.  In one
example we will be able to prove the analogs of these conjectures.

This paper is structured as follows.
In Section 2,
we first recall some facts about white noise measures.  These
are the infinite dimensional analogs of the Gaussian measures
$\mu_W$ above.
We then study maps $\cF: X \to Y$ where $X$ is a
Banach manifold and $Y$ is a Banach space
equipped with a white noise measure.  If
$\cF :X \to Y$ is a differentiable (that is, $C^1$) map,
we obtain a measure on $X$
analogous to the measure $f^*\mu_W$ arising in (\ref{measdef}).
This is done in Proposition \ref{main}. 

In Section 3,
we give examples of this construction of measures where the maps are motivated
by supersymmetric quantum field theories appearing in the physics
literature.  In these examples we also give a construction of the phase of
the determinant using a generalization of Fredholm determinants.
We study three families of measures.  The first family is denoted
$\mu_P^{s,t},$ where $P$ is a polynomial in one variable and $s,t$ are
parameters.  These are measures on a Banach space of functions
on the torus, and are
related to objects known as the Wess-Zumino-Landau-Ginzburg
models (Theorem \ref{wzlgthm}).  We next consider a family
$\mu_\cG^{s,t}$ of measures on a group $\cG$ of maps from $\P^1$ to a complex
Lie group; these measures
are related to the Wess-Zumino-Novikov-Witten models (Theorem 
\ref{wznwthm}).  Finally, given a compact Lie group $G$ and a three-manifold
$M,$
we consider measures $\mu_{M,G}^{s,t}$
on the product of a space of connections on a principal
$G$ bundle on $M$ with a space of $\fg$-valued three-forms on $M.$
These measures are related to three-dimensional quantum gauge theory
(Theorem \ref{3dgtthm}).
We conjecture that the measures $\mu_P^{s,t}$ and $\mu_{M,G}^{s,t}$
extend to Borel probability measures.  We also give precise conjectures about
the values of the integrals of the phases of appropriate determinants
(Conjecture \ref{wzlgconj} and Conjecture \ref{casson}).  In the
case of three-dimensional gauge theory the conjecture is that
$\mu_{M,SU(2)}^{s,t}$ should yield the Casson invariant of homology
three-spheres.

At the end of this section we briefly describe how our methods
can be used in a few other examples of maps of Banach manifolds
inspired by supersymmetric quantum field theories.

In Section 4 we turn to the measures $\mu_\cG^{s,t},$
where we are able to prove analogs of the conjectures we
have made about $\mu_P^{s,t}$ and $\mu_{M,G}^{s,t}.$
Using results of Atiyah and Bott \cite{ab}, we show that $\mu_\cG^{s,t}$ is a Borel
probability measure on the space $\cG$ (Theorem \ref{emb}).
We then compute (Theorem \ref{fzt}) some
expectations of functions in this measure to recover
formulas discovered by Frenkel and Zhu \cite{fz} in the context of
vertex operator algebras.

\section{Maps and measures on Banach manifolds}

In this section we consider some general constructions of measures
on Banach manifolds.  We begin with Gaussian measures and then present
a method of producing non-Gaussian measures from Gaussian measures
and maps of Banach manifolds.\footnote{In this paper we consider both the case of complex Banach spaces and
that of real Banach spaces.  The results presented in this section
hold in either the real or complex category, and all formulas should
be interpreted in the appropriate setting.}

\subsection{White Noise and Gaussian Measures}

Let $M^d$ be a smooth compact Riemannian manifold and let $V \to M$
be a smooth Hermitian vector bundle over $M.$
Let $\Gamma(V)$ denote the space of 
smooth sections of $V.$  The Hermitian metric on $V$ and
the Riemannian structure on $M$ give rise to an $L_2$ inner product
$<,>$ on $\Gamma(V)$.  For $\l \in \R,$ let $H_\l(M)$
denote the distributions on $M$ of Sobolev class $\l$ and let 
$\G_\l(V) = \G(V) \otimes H_\l(M)$ denote the space of sections
of $V$ of Sobolev class $\l.$  Let $\alpha$ be a connection on $V,$
and let $\Delta_\alpha$ be the corresponding Laplacian on sections
of $V.$  This Laplacian gives rise to
a hermitian inner product on $\G(V)$ given by
$\langle s,s'\rangle_\l = <s , (-\Delta_\alpha + 1)^{\lambda}s'>$
for $s,s' \in \G(V).$  This inner product extends to $\G_\l(V)$
and makes $\G_\l(V)$ into a separable Hilbert space.
If $\l > 0$ the
space $\G_{-\l}(V)$ may be identified
with the dual of $\G_\l(V);$ we continue to denote the
pairing of $\G_{-\l}(V)$ with $\G_\l(V)$ by $<,>.$
If $\l=0$ we write $\G_0(V) = L_2(M,V),$ and the pairing
$<,>$ gives rise to the $L_2$ metric on $\G_0(V).$

If $\lambda < -d/2,$ the space $\G_\l(V)$
comes equipped with a finite Borel measure $\mu_V$ called
{\em white noise measure,} which
is characterized by the following property:

Let $\phi \in \G_{-\l}(V)$ and let $E_\phi : \G_\l(V) \to \C$ denote the
function given by 
$$ E_\phi(\sigma) = \exp(i{\rm Re} <\phi,\sigma>).$$

\noindent  Then 
\begin{equation}\labell{char}
 \mu_V(E_\phi) = \exp( -\frac12 ||\phi||^2_{L_2}).\end{equation}

In particular, $\mu_V(1) = 1.$

Given $t > 0,$ let $R_t : \G_\l(V) \to \G_\l(V)$ be the scaling
function $R_t(v) = t v.$  We write $\mu_V^t = (R_{t^{-1}})_* \mu_V$ for
the scaled white noise measure.  Then for every $\epsilon > 0,$ 
we have
$\mu_V^t(B_\epsilon(0)) > 0$ for $t$ sufficiently large.

We now define more general Gaussian measures. Choose
$s\in \R.$  Then $(-\Delta_\alpha + 1)^{-s}:
\G_\lambda(V)\to \G_{\lambda+ 2s}(V) $ is an isomorphism.  Let $\lambda
< -d/2 $ as before; then we may push the white noise measure $\mu_V^t$
forward to obtain a Gaussian measure
$\mu_{V,s}^t:= ((-\Delta_\alpha + 1)^{-s})_* \mu_V^t$ on $\G_\nu(V)$ where
$\nu < 2s - d/2.$\footnote{
In the literature these measures are usually regarded
as living on the spaces $Z_\lambda(V)$ of sections of $V$ of
{\em Zygmund} class $\l < 2s-d/2.$  However, an extension
of the Sobolev embedding theorem (see \cite{hormander}, Proposition
8.6.10) shows that there exist continuous inclusions
$$Z_{\lambda + \epsilon}(V) \subset \G_\l(V) \subset Z_{\l- d/2}(V)$$
\noindent for any $\epsilon > 0,$ so that our measures can just as
well be regarded as living on $\G_\l(V).$  This will be convenient
for various purposes in this paper.
}

\begin{Proposition}[Cameron-Martin]\labell{CM}Let $\lambda < 2s -d/2$ and
let $v \in \G(V).$  Let $T_v:\G_\l(V) \to \G_\l(V)$ be the translation
given by $T_v(a) = v + a.$  Then $({T_v})_*\mu_{V,s}^t$ is absolutely
continuous with respect to $\mu_{V,s}^t,$ and the Radon Nikodym
derivative may be written as

\begin{equation}
\frac{d (T_v)_*\mu_{V,s}^t}{d\mu_{V,s}^t} = 
\exp(t^2<(-\Delta_\alpha + 1)^{2s} v,\cdot> -\frac12 t^2 ||v||^2_{2s}).
\end{equation}
\end{Proposition}
 
The translate of the measure $\mu_{V,s}^t$ by an arbitrary 
element $v\in \G_\l(V)$
may not be absolutely continuous with respect to $\mu_{V,s}^t.$   

\subsection{Pullbacks of measures by differentiable maps and the main
construction}\labell{tmc}

Let $X$ be a separable metric space.  
Let $\{ (U_\a , \mu_\a)\}_{\a \in \cA}$ be a collection
of pairs $(U_\a, \mu_a)$ where $\{U_\a\}_{\a\in \cA}$ is an open cover
of $X$ and $\mu_\a$ is a finite Borel measure on $U_\a$ for each $\a.$
Suppose that for each $\a,\beta \in \cA,$ 

$$ \mu_\a |_{U_\a \cap U_\beta} = \mu_\beta|_{U_\a \cap U_\beta} .$$  

\noindent Since $X$ is Lindel\"of, we obtain a $\sigma$-finite measure $\mu$
on
$X.$\footnote{To do this, extend the measure to any countable subcover
of the cover $\{U_\alpha\}_{\alpha \in \cA}.$ This extension is independent
of the choice of countable subcover.}
We will call such a measure a {\em Borel
probability measure} if $\infty > \mu(X) > 0.$

Let $X,Y$ be separable Banach manifolds
and let $\mu$ be a finite Borel measure on $Y.$
Suppose $\cF : X \to Y$ is a differentiable
map, and let $U\subset X$ be the (open) set of points $p \in X$
such that $\delta \cF_p$ is an isomorphism.  For each $p\in U$,
we can find a neighborhood $U_p$ of $p$ where $\cF|_{U_p}$ is
invertible.  Denote this local inverse by $\cG_p,$  and define a finite Borel
measure $\nu_{U_p}$ on $U_p$ by
$$\nu_{U_p} := {\cG_p}_* (\mu|_{\cF(U_p)}).$$
\noindent Then the collection $\{(U_p,\nu_{U_p})\}$ gives a measure
$\nu$ on $U.$  Let $i : U \to X$ be the inclusion.
We define the measure $\cF^*\mu$ on $X$ by 
$\cF^*\mu  = i_* \nu.$

We therefore obtain the following result, which
is our main tool for constructing
measures on Banach manifolds.

Recall that a differentiable map $\cF : X \to Y$ of separable Banach spaces is
a {\em Fredholm map} if $\delta \cF|_x$ is Fredholm for all $x \in X.$

\begin{Proposition}\labell{main}
Let $X,Y$ be separable Banach manifolds.  Let 
$\mu$ be a finite Borel measure on $Y.$  Let $\cF:X\to Y$
be a differentiable map.
Then the map $\cF$ induces a $\sigma$-finite Borel measure $\cF^*\mu$
on $X,$ supported on the open set $U \subset X$ where $\delta \cF$ is
an isomorphism.

In particular, let $M^d$ be a compact
Riemannian manifold of dimension $d,$ and let $V$ be a vector
bundle over $M.$  Let $\l < -d/2,$ and
let $\cF: X \to \G_\l(V)$ be a differentiable map.
Suppose
that $\cF$ is a Fredholm map and that ${\rm ind~}\delta \cF|_x = 0$
for all $x \in X.$
Then for any $t > 0$ the map $\cF$ induces a $\sigma$-finite Borel measure
$\cF^*\mu_V^t$ on $X,$ supported on the set of regular points of $\cF.$

Suppose that there exists $p \in X$  which is a regular point of $\cF$ and
satisfies $\cF(p) = 0.$  Then for $t$ sufficiently large,
there exists a neighborhood $W \subset X$
of $p$ such that $$\infty > \cF^*\mu_V^t(W) > 0.$$
\end{Proposition}

The methods we have used in this section cannot decide whether $\cF^*\mu_V^t$
gives a Borel probability measure on $X.$  This has to be
checked by analysis in each example.  We will do this in one example
in Section 4.

\section{Examples of measures on Banach manifolds}
In this section we give examples of maps which can
be used to construct measures
on Banach manifolds by the construction of
Proposition \ref{main}.  These maps are inspired
by formal path integrals appearing in the physics literature.

\subsection{Some analytical preliminaries}
In this section we collect some analytical facts we will use
in our constructions.

\subsubsection{Sobolev bounds}
In order to study nonlinear maps on Banach manifolds,
we will use the following
result on the regularity of a product of distributions.

\begin{Proposition}\labell{sob}[see e.g. \cite{hormander}, Prop. 8.3.1]
\footnote{ We will use the bound (\ref{bound}) only in the
situation where $s_1, s_2 > d/2;$ in other words,
when the Sobolev spaces $H_{s_i}(M^d)$ are spaces of
continuous functions. }

Let $f_1 \in H_{s_1}(M^d),$ $f_2 \in H_{s_2}(M^d).$  Suppose
$s_1+ s_2 \geq  0.$  Then $f_1 f_2 \in H_s(M^d)$ if 

$$s \leq s_i, ~i = 1,2 {\rm ~and ~} s < s_1 + s_2 -d/2.$$

\noindent The product $f_1 f_2$
is bounded by 
\begin{equation}\labell{bound}
|| f_1 f_2 ||_s <  K ||f_1||_{s_1} ||f_2||_{s_2}.
\end{equation}
\noindent where $K$ is a constant depending on $s,s_1,s_2.$
\end{Proposition}

This result immediately implies a similar result for sections
of vector bundles of the appropriate Sobolev classes.

\subsubsection{Schatten classes and regularized determinants}

To define phases of determinants of some of the operators we
will encounter, we will use regularized determinants \cite{ds,simon}.

Let  ${\mathcal H}$ be a separable Hilbert space, and let
$\cL(\cH)$ denote the ring of bounded operators on $\cH.$  For $k \geq 1,$
write 
$\cS_k({\mathcal H})$ for the $k$-th Schatten ideal of ${\mathcal H};$
that is, the set of compact operators $K:{\mathcal H}\to {\mathcal H}$
such that $(K^*K)^{k/2}$ is trace class.
The space $\cS_k(\cH)$ is a two-sided ideal of $\cL(\cH).$
The function
$K\to ||K||_k:={\rm tr}(K^*K)^{k/2} $ defines a norm on
$\cS_k({\mathcal H})$ which makes $\cS_k({\mathcal H})$ into a
Banach space, and the inclusion $\cS_k(\cH) \to \cL(\cH)$ is continuous.
Similarly, if $k < k',$ $\cS_k(\cH) \subset \cS_{k'}(\cH),$
and the inclusion $\cS_k(\cH) \to  \cS_{k'}(\cH)$ is continuous.

An example of an element of a Schatten ideal is given as follows.
Let $M$ be a compact Riemannian manifold of dimension $d$ and
let $V$ be a Hermitian vector bundle over $M.$  Let $\alpha$ be a connection
on $V,$ let
$(-\Delta_\alpha + 1): \G_{\l+2}(V) \to \G_{\l}(V)$ be the Laplacian, and
let $j: \G_{\l+ 2}(V) \to \G_\l(V)$ denote the inclusion.
Then
\begin{equation}\labell{laplace}
j \circ (-\Delta_\alpha + 1)^{-1}\in \cS_{r} (\G_{\l}(V))
\end{equation}
\noindent for all $r > d/2.$

Recall that if $K \in \cS_1(\cH),$ the {\em Fredholm determinant}
$det(1 + K)$ is defined  \cite{groth} by the convergent series
\begin{equation}\labell{fredholm}
det(1+K) := \sum_{n=0}^\infty {\rm tr~} \wedge^n K.
\end{equation}
This construction was generalized to operators in $\cS_k(\cH)$ by
Poincar\'e \cite{poincare}; we follow Simon \cite{simon}.  Given
$K \in \cL({\cH}),$ define $R_k(K) \in \cL(\cH)$ by
\begin{equation}\labell{rkk}
R_k(K) := [(1 + K) \exp (\sum_{n=1}^{k-1} \frac{(-K)^n}{n})] -1.
\end{equation}
Simon shows that if $K \in \cS_k(\cH),$ then $R_k(K) \in \cS_1(\cH).$
Define the {\em regularized determinant}  $det_k(1+K)$ by\footnote{The motivation for this definition is  Plemelj's formula  (valid when $K \in \cS_1(\cH)$ and $||K||_1 < 1$)
$$det (1+K) = \exp ({\rm tr~} \log(1+ K)).$$
The appearance of $R_k(K)$ in (\ref{regdet}) amounts to
removing the first $k-1$ terms from the power series
for $\log (1+K)$ in Plemelj's formula.
The remaining terms  are multiples of $K^n$ where
$n\geq k$ and so are trace class if $K \in S_k(\cH).$}
\begin{equation}\labell{regdet}
det_k(1+K) := det(1 + R_k(K)).
\end{equation}
The function
$det_k(1 + \cdot) : \cS_k({\mathcal H})\to \C$ is continuous.
If $K$ is trace class,
the regularized determinant is related to the Fredholm determinant
$det(1 + K)$ by 
\begin{equation}\labell{trcl}
 det_k(1 + K) = det(1 + K) \exp ( \sum_{n=1}^{k-1} \frac{(-1)^n}{n}{\rm tr~ } K^n).
\end{equation}

The regularized determinant has multiplicativity properties which
can be deduced from (\ref{trcl}) and the multiplicativity
of the Fredholm determinant by approximating elements of the
Schatten classes by trace-class operators (which are dense
in $\cS_k({\mathcal H})$).  From these multiplicativity
properties it follows that if $1 + K$ is invertible, the
regularized determinant $det_k(1+K)$ is nonzero.

\subsection{Wess-Zumino-Landau-Ginzburg models on the torus}  
Let $\Sigma = S^1 \times S^1$ be the 
torus, equipped with a K\"ahler structure, and let $V =  \Sigma\times \C \to \Sigma$
be the trivial bundle over $\Sigma.$  If $\l \in \R,$ the space
$\G_\l(V)$ is given by  $\G_\l(V)= H_\l(\Sigma).$  Choose $\l < -1$
and $s > -\l/2.$  Let $\mu_V$ be the white noise measure on $\G_\l(V).$
By the Sobolev embedding theorem, the space
$\G_{\l + 2s +1} (V)$ consists of continuous functions.

Let $P$ be a polynomial in one variable.  Define 
$\cF_s:\G_{\lambda+2 s+ 1}(V) \to \G_\lambda(V)$ 
by

\begin{equation}\labell{wzlgs}
\cF_s(\phi) = (-\Delta + 1)^s({\partial} \phi + \overline{P'(\phi)}).
\end{equation}

\begin{Theorem}\labell{wzlgthm}
The map $\cF_s$ is well-defined
and differentiable. Its derivative is a
Fredholm operator of index $0.$

Therefore, for each $t > 0,$ the pullback measure
$$\mu_P^{s,t}:= \cF_s^*\mu_V^t$$
exists and is a $\sigma$-finite measure on $\G_{\l+ 2s +1}(V),$
supported on the set of regular points of $\cF_s.$
If the polynomial $P'$ has a transversal zero at $y \in \C,$ then for
$t$ sufficiently large, there exists a neighborhood $W$
of the constant function $\phi_0 = y$
such that $$\infty > \mu_P^{s,t}(W) > 0.$$
\end{Theorem}

\begin{proof}  Proposition \ref{sob} shows that the
map $\cF_s$ is well-defined and differentiable.
The derivative $\delta \cF_s$ is a compact perturbation
of the isomorphism
$(-\Delta+1)^s({\partial} + 1):\G_{\lambda+2 s+ 1}(V) \to \G_\lambda(V),$
and so is Fredholm of index $0.$  The existence of the measure
follows by Proposition \ref{main}.

If $y\in \C$ is a point where $P'(y) = 0$ and
$P''(y) \neq 0,$ the constant function $\phi_0 :=y$ satisfies
$\cF_s(\phi_0) = 0,$ and $\delta \cF_s|_{\phi_0} = 
(-\Delta+1)^s({\partial} + P''(y))$ is an isomorphism.  The
positivity of the measure $\cF_s^*\mu_V^t$ in a neighborhood of
$\phi_0$ for $t$ sufficiently large follows, again, by Proposition \ref{main}.
\end{proof}

We make the following conjecture.

\begin{Conjecture}\labell{convwzlg}
The measure $\mu_P^{s,t}$ satisfies  $$ \mu_P^{s,t}(1) < \infty$$
for $t$ sufficiently large.
\end{Conjecture}

Suppose the polynomial
$P'$ has a transversal zero. Then for
$t$ sufficiently large, Conjecture \ref{convwzlg} implies that 
the measure $\mu_P^{s,t}$ is a Borel probability
measure on $\G_{\l + 2s  +1}(V).$

\subsubsection{The phase of the determinant of $\delta \cF_s$}
The derivative $\delta \cF_s|_{\phi}$ is given by the Fredholm operator
$$D_\phi:=(-\Delta+1)^s({\partial} + \overline{P''(\phi)}):
\G_{\lambda+2 s+ 1}(V) \to \G_\lambda(V).$$ 
If $P'$ is nonconstant
and $y$ is a point where $P''(y)\neq 0,$ then if $\phi_0:=y$ is
the constant function, $D_{\phi_0}$ is an isomorphism.  Then
$(D_{\phi_0})^{-1} D_\phi :\G_{\l+ 2s +1}(V) \to\G_{\l+ 2s +1}(V)$
can be written as a compact perturbation of the identity 
$$(D_{\phi_0})^{-1}D_\phi = 1 + K(\phi_0,\phi),$$ where
$$K(\phi_0,\phi) = 
(\partial + P''(\phi_0))^{-1}\circ j \circ M(\overline{P''(\phi)} - \overline{P''(\phi_0)}).$$

\noindent Here for $\psi \in \G_{\l + 2s + 1}(V),$ we have denoted by 
$$M(\psi): \G_{\l + 2s + 1}(V)\to \G_{\l + 2s + 1}(V)$$ the operator
given by multiplication by $\psi,$  
$j: \G_{\l + 2s + 1}(V)\to \G_{\l + 2s}(V)$ is the inclusion, and
$(\partial + P''(\phi_0))^{-1}$ is the inverse of the isomorphism
$$\partial + P''(\phi_0): \G_{\l + 2s + 1}(V)\to \G_{\l + 2s}(V).$$
Since $||M(\psi)|| \leq ||\psi||_{2s +\l +1}$ by Proposition \ref{sob},
the estimate (\ref{laplace}) shows that
$$K(\phi_0,\phi) \in \cS_3(\G_{\l+ 2s +1}(V))$$
for all
$\phi.$  The regularized determinant $det_3(1+K(\phi_0,\phi))$ is therefore
well-defined for all $\phi,$ and is nonzero at all regular points
of $\cF_s.$  Let $U$ be the set of regular points of $\cF_s.$
Define $\Psi(\phi_0,\cdot): U \to S^1$ by
\begin{equation}\labell{detwzlg}
\Psi(\phi_0,\phi)=\frac{det_3(1+K(\phi_0,\phi))}{|det_3(1+K(\phi_0,\phi))|}.
\end{equation}
Morally, the function $\Psi(\phi_0,\phi)$
gives the ratio of the ``phase of the determinant
of $\delta\cF_s|_\phi$'' to the
``phase of the determinant of $\delta\cF_s|_{\phi_0}.$''
This method of defining determinants is used in the physics
literature;
by analogy with the case of maps of finite-dimensional vector spaces
given in equations (\ref{ls}) and (\ref{phasofmeas}) in
the introduction, we make the following conjecture.
\begin{Conjecture}\labell{wzlgconj}
Assume that Conjecture \ref{convwzlg} holds.  Then for $t$ sufficiently
large,
$$|\mu_P^{s,t} (\Psi(\phi_0,\cdot) )| =  ({\rm deg~} P - 1).$$
\end{Conjecture}

Morally this conjecture says that for $t$ large, the measure
$\mu_P^{s,t}$ is concentrated at the
solutions of the nonlinear equation $\cF_s = 0.$

\begin{Remark}The formal path integral for
the Wess-Zumino-Landau-Ginzburg model appearing in the physics literature
corresponds to formula (\ref{pb}) applied (formally) to $\cF_0,$ and should
morally count the solutions of the partial differential
equation  $\cF_0 = 0.$  Since the solutions
of the equation $\cF_s = 0$ are exactly those of the equation 
$\cF_0 = 0,$ the computations associated formally
with the path integral corresponding to $\cF_0$ should morally
be replicated for the measure $\mu_P^{s,t}.$
The supersymmetric field theory corresponding
to these path integrals was constructed using classical
methods of constructive quantum field theory in \cite{jlw},
and the resulting partition function does count (for 
generic $P$) the number of solutions of $\cF_0 = 0,$ which correspond
to constant functions with values given by the zeros
of $P.$

The idea of using an infinite-dimensional
degree to find solutions of nonlinear partial differential
equations goes back to Leray-Schauder, who considered {\em
proper} maps.  The case of measures induced by nonlinear maps on
$\G_\l(V)$ where $V$ is a vector bundle over a one-dimensional
manifold was considered in \cite{getzler1} using the methods of 
\cite{kuo,ramer,kusuoka,getzler}.
In this case the pullback measure on $\G_\l(V)$ is a perturbation
of a Gaussian measure, and the analog of our operator
$K$ is Hilbert-Schmidt.   This is not the case if the
underlying manifold is of dimension greater than one.

The subtlety in the definition of the phase of the determinant
arising from the need to use regularized determinants was to my
knowledge first observed in the physics literature.\footnote{The appearance of regularized determinants is one reason
we have chosen to work with Sobolev spaces and Hilbert manifolds
rather than with Zygmund spaces and Banach manifolds.  While there
exist theories of determinants on Banach spaces (see e.g. \cite{groth,ruston})
I am not aware of a theory of regularized determinants in this
context.}
\end{Remark}

\subsection{Wess-Zumino-Novikov-Witten models on $\P^1$}

Let $G$ be a compact semisimple Lie group.  We consider
the corresponding complex Lie group as a subgroup
$G^\C\subset SL(n,\C)$ of $SL(n,\C).$
Then $G^\C$ and $\fg^\C \subset {\mathfrak sl}(n,\C)$
are subsets of the space
$M_n(\C)$ of $n \times n$ matrices, and $\fg^\C$ acquires
a Hermitian inner product from this inclusion.\footnote{The only
reason for the restriction $G^\C \subset SL(n,\C)$ is to allow products
and
the inverse map on the space $\cG$ to be estimated directly using 
Proposition \ref{sob}.
It is not very difficult to remove this restriction; see \cite{palais}.}
Consider $\P^1,$ equipped with the standard
K\"ahler metric, and choose a base point in $\P^1.$  Let $\lambda < -2$
and $ s > (-\lambda+ 1)/2, $  and let
$\cG:={\rm Map}_{*,\lambda+ 2s+ 1}(\P^1, G^\C)$
be the space of maps from $\P^1$
to $G^\C$ of Sobolev class $\lambda + 2s+ 1$ which carry the base point
in $\P^1$ to the identity in $G^\C.$   We have
chosen $s$ sufficiently large so that the space $\cG$ is a Lie group, with
the group structure given by pointwise multiplication.  
As a Hilbert manifold $\cG$ is modelled on the Lie algebra ${\rm Lie}(\cG)$
given by the space
${\rm Lie}(\cG) :
= \Omega^{0}_{*,\lambda + 2s + 1} (\P^1, \fg^\C):=\G_{*,\l +2s + 1}(V)$ of
sections of the trivial bundle $V := \P^1 \times \fg^\C \to \P^1$ 
of Sobolev class $\l+ 2s +1$ which vanish at the base point of $\P^1.$
For $\nu \in \R,$ let $\Omega^{0,1}_\nu(\P^1, \fg^\C) := 
\Omega^{0,1}(\P^1, \fg^\C) \otimes H_\nu(\P^1)$
be the space of $\fg^\C$-valued
$(0,1)$-forms on $\P^1$ of Sobolev class $\nu;$
if $W =  {T^*\P^1}\otimes\fg^\C,$
then $\Omega^{0,1}_\nu(\P^1, \fg^\C)  = \G_\nu(W).$
The space
$\Omega_{\lambda}^{0,1}(\P^1,\fg^\C) = \G_\l(W)$ 
is equipped with a white noise measure $\mu_W.$

The group $\cG$ acts freely on the space $\cA_{\l + 2s} 
:=\Omega^{0,1}_{\l + 2s}(\P^1, \fg^\C)$ by the formula
$g \cdot A = g A g^{-1} + \bar{\partial}g g^{-1}$ for
$g \in\cG,$ $A \in \cA_{\l + 2s}.$
Let $\Delta : \Omega^{0,1}_{\nu + 2}(\P^1, \fg^\C)\to\Omega^{0,1}_\nu(\P^1, \fg^\C)$ denote the Laplacian.

Define the map
$\cF_s : \cG \to \Omega^{0,1}_{\lambda} (\P^1, \fg^\C)$
by
\begin{equation}\labell{wznws}
\cF_s(g) = (-\Delta + 1)^s \bar{\partial} g g^{-1}
\end{equation}
\begin{Theorem}\labell{wznwthm}
The map $\cF_s$ is well-defined
and  differentiable. Its derivative is a
Fredholm operator of index $0.$  The map
$\delta\cF_s|_g$ is an isomorphism for all 
$g \in \cG.$

Therefore, for each $t > 0,$ the pullback measure
$$\mu_\cG^{s,t}:= \cF_s^*\mu_W^t$$
exists and is a $\sigma$-finite measure on $\cG.$
 For $t$ sufficiently large, this measure is positive.
\end{Theorem}

In the next section we will prove that the measure $\mu_\cG^{s,t}$ is
a Borel probability measure for all $t;$ see Theorem \ref{emb}.

\begin{proof}  
An application of Proposition \ref{sob} shows that $\cF_s$ is
well-defined and differentiable.  For $g \in \cG,$ the tangent
space $T\cG_g$
can be identified with
$\Omega^{0}_{*,\lambda+ 2s + 1} (\P^1, \fg^\C).$
The map
$\cF_s$ is the composite of the map $\cF: \cG \to \cA_{\l+ 2s}$
given by $\cF(g) = \bar{\partial}g g^{-1}$ with the
isomorphism $(-\Delta + 1)^s:\cA_{\l + 2s} \to
\Omega^{0,1}_{\lambda} (\P^1, \fg^\C).$
The map $\cF$ gives a smooth
action of $\cG$ on $\cF(\cG),$ which is the restriction to $\cF(\cG)$
of the action of $\cG$ on $\cA_{\l + 2s}.$ 
In terms of this group action,
$\cF(\cG)$ is the orbit of the point $0 \in \cA_{\l + 2s}.$
The derivative of $\delta \cF|_e$ at
the identity $e \in \cG$ is $\delta\cF|_e = \bar{\partial} : 
\Omega^{0}_{*,\lambda+ 2s + 1} (\P^1, \fg^\C) \to 
\Omega^{0,1}_{\lambda + 2s} (\P^1, \fg^\C),$ which is an isomorphism.
Therefore the derivative of $\delta\cF|_g$ is an isomorphism at any point
$g \in \cG.$ Compare \cite{ab}, Lemma 14.6 and (14.7).
\end{proof}

\subsubsection{The phase of the determinant of $\delta \cF_s$}

The derivative $\delta \cF_s|_{g}$ is given by the Fredholm operator
$$D_g:=(-\Delta+1)^s(\bar{\partial} - \bar{\partial}g g^{-1}):
\Omega^{0}_{*,\lambda+ 2s + 1} (\P^1, \fg^\C)
\to
\Omega^{0,1}_{\lambda} (\P^1, \fg^\C).$$  Choose a smooth base point
$g_0 \in \cG.$  Then
$(D_{g_0})^{-1}D_g :
\Omega^{0}_{*,\lambda+ 2s + 1} (\P^1, \fg^\C)
\to
\Omega^{0}_{*,\lambda+ 2s + 1} (\P^1, \fg^\C)$
can be written as a perturbation of the identity 
$$(D_{g_0})^{-1}D_g = 1 + K(g_0,g),$$
 where
$$K(g_0,g) = 
(\bar{\partial} - \bar{\partial}g_0 g_0^{-1})^{-1}\circ
M(\bar{\partial}g_0g_0^{-1} -
\bar{\partial}g g^{-1}).$$

\noindent Here for $\psi \in \Omega^{0,1}_{\l + 2s}(\P^1, \fg^\C),$ we
have denoted by $$M(\psi): \Omega^{0}_{*,\lambda+ 2s + 1} (\P^1, \fg^\C)
\to \Omega^{0,1}_{\lambda+ 2s} (\P^1, \fg^\C)$$
the wedge product by $\psi,$ and $(\bar{\partial} - \bar{\partial}g_0 g_0^{-1})^{-1}$
is the inverse of the isomorphism
$$(\bar{\partial} - \bar{\partial}g_0 g_0^{-1}):\Omega^{0}_{*,\lambda+ 2s + 1} (\P^1, \fg^\C)
\to \Omega^{0,1}_{\lambda+ 2s} (\P^1, \fg^\C).$$

We will show that
$K(g_0,\cdot) \in \cS_3(\Omega^{0}_{*,\l + 2s+1} (\P^1, \fg^\C))$
almost everywhere with respect to the measure $\mu_\cG^{s,t}.$
Since $\l < -2,$ the white noise measure on $\G_\l(W)$ is supported
on $\G_{\l + 1}(W) \subset \G_\l(W).$  Thus
the pullback measure $\mu_\cG^{s,t}=\cF_s^*\mu_W^t$ is supported
on ${\rm Map}_{*,\l + 2s + 2}(\P^1,G^\C)\subset \cG.$  It follows that 
$(\bar{\partial}g_0g_0^{-1} - \bar{\partial}g g^{-1}) \in 
\Omega^{0,1}_{\l + 2s +1}(\P^1,\fg^\C)$ for $g$ lying in the complement
of a set of measure zero in the measure $\mu_\cG^{s,t}.$

Then, for almost all $g,$
$$K(g_0,g) = 
(\bar{\partial} - \bar{\partial}g_0 g_0^{-1})^{-1}\circ j \circ
\hat{M}(\bar{\partial}g_0g_0^{-1} -
\bar{\partial}g g^{-1}).$$

\noindent where for $\psi \in \Omega^{0,1}_{\l + 2s + 1}(\P^1, \fg^\C),$ we
have denoted by $${\hat{M}}(\psi): \Omega^{0}_{*,\l + 2s +1} (\P^1, \fg^\C)
\to \Omega^{0,1}_{\l + 2s +1} (\P^1, \fg^\C)$$
the wedge product by $\psi,$ 
and $j: \Omega^{0,1}_{\l + 2s + 1} (\P^1, \fg^\C) \to \Omega^{0,1}_{\l + 2s} (\P^1, \fg^\C)$
is the inclusion.

Then by Proposition \ref{sob}, $||{\hat{M}}(\psi)|| \leq ||\psi||_{\l + 2s+1},$
and by the estimate (\ref{laplace}),
$${K}(g_0,g) \in \cS_3(\Omega^{0}_{*,\l + 2s+1} (\P^1, \fg^\C))$$
for almost all $g.$
The regularized determinant $det_3(1+K(g_0,g))$ is therefore
well-defined almost everywhere, and since $1+K(g_0,g)$ is invertible,
$det_3(1+K(g_0,g))$ is nonzero for all $g$ where it is defined. Define 
$\Psi(g_0,\cdot): \cG \to S^1$ almost everywhere by
\begin{equation}\labell{phasewznw}
\Psi (g_0,g) = \frac{det_3(1+K(g_0,g))}{|det_3(1+K(g_0,g))|}.
\end{equation}

\begin{Remark}\labell{heuristic}  The expression (\ref{trcl}) for the determinant
of a trace-class operator gives a conjecture for the value
of the function $\Psi.$  Approximate $K$ by a sequence $K_n$
of trace-class operators.  
Since $\cF_s$ is a holomorphic map,
we expect that morally the phase of the determinant $det(1 + K_n)$
will tend to $+1$ as $n \to \infty.$
In view of equation (\ref{trcl}) the phase of the regularized determinant $det_3(1+K)$
should be given morally by
$\Psi = \exp(i \lim_{n\to \infty}{\rm Im~tr }( {K_n}^2)),$
so that we might expect $\Psi$ to be the exponential of
a (possibly nonlocal) quadratic polynomial in $\bar{\partial}g g^{-1}.$
This quantity is known as the anomaly
in the physics literature, and often turns out to
be well-defined when $K$ is taken as a limit of trace-class
operators, even though $K$ is not Hilbert-Schmidt.\footnote{The
physics literature uses a type of approximation of $K$
by a trace-class operator called dimensional regularization.}
See also Remark \ref{phirmk}.
\end{Remark}
\subsection{Three-dimensional Gauge Theory}  
Let $M$ be a compact, smooth Riemannian 3-manifold.
Let $G$ be a compact Lie group, and choose an invariant Hermitian
inner product on $\fg.$  Choose $\lambda < -3/2$  and $s >  (-\lambda+1/2)/2$
and let $\cA= \cA_{\lambda+2s+1}$ be the
space of connections on the trivialized principal $G$-bundle on $M$
of Sobolev class $\lambda+2s +1;$ in the notation of Section 2,
this is given by $\cA = \G_{\l+2s +1}(T^*M\otimes \fg).$
For $A \in \cA$ denote the curvature of $A$ by $F_A.$
After a choice of a base point in $M,$ connections
$A$ on $M$ which satisfy the equation $F_A = 0$ correspond 
by the monodromy representation to
representations of $\pi_1(M)$ in $G.$
Choose such a flat connection $A_0.$
Let $\Omega_{\nu}^i(M,\fg)= \G_\nu(\Lambda^i(T^*M)\otimes \fg)$
be the space of de Rham forms
on $M$ with values in $\fg$ of Sobolev class $\nu.$
The space
$\Omega_{\lambda}^0(M,\fg)\oplus\Omega_{\lambda}^2(M,\fg)$ is given,
in the notation of Section 2, by $\G_\l(V)$ where
$V = \fg \oplus (\Lambda^2 T^*M)\otimes \fg,$ and where $\l < -3/2.$
It is therefore equipped with a white noise measure $\mu_V.$

Denote by $*$ the Hodge $*$ operator on differential forms
with values in
$\fg$ obtained from the Riemannian pairing on differential forms
on $M$ and the hermitian
inner product on $\fg.$  
Consider
the map $\cF_{s,A_0} : \cA_{\lambda+2s+1} \oplus \Omega_{\lambda+2s+1}^3(M,\fg)
\to \Omega_{\lambda}^0(M,\fg)\oplus\Omega_{\lambda}^2(M,\fg)$ be given by

\begin{equation}\labell{3dgauge}
\cF_{s,A_0}(A,\xi) = (-\Delta_{A_0} + 1)^s (*d_{A_0}*(A-A_0) + F_A + *d_A*\xi).
\end{equation}

\begin{Theorem}\labell{3dgtthm}
The map $\cF_{s,A_0}$ is well-defined
and differentiable. Its derivative is a
Fredholm operator of index $0.$
Therefore, for each $t > 0,$
the pullback measure
$$\mu_{M,G}^{s,t,A_0}:=\cF_{s,A_0}^*\mu_V^t$$
exists and is a $\sigma$-finite measure on
$ \cA_{\lambda+2s+1} \oplus \Omega_{\lambda+2s+1}^3(M,\fg),$
supported on the set of regular points of $\cF_{s,A_0}.$

Suppose that $A_0$ is chosen to be a flat connection that corresponds
to an isolated conjugacy class of
irreducible representations of $\pi_1(M)$ in $G.$
Then for $t$ sufficiently
large, there exists a neighborhood $W$ of $(A_0,0)$
such that $$\infty > \mu_{M,G}^{s,t,A_0}(W) > 0.$$
\end{Theorem}

\begin{proof}  Again, we can use Proposition \ref{sob} to show that
$\cF_{s,A_0}$ is well-defined and differentiable.
The derivative $\delta \cF_{s,A_0}$ is a compact perturbation
of the operator
$ (-\Delta_{A_0}+1)^s(*d_{A_0}* + d_{A_0}):
\Omega_{\lambda+2s+1}^1(M,\fg) \oplus \Omega_{\lambda+2s+1}^3(M,\fg)
\to \Omega_{\lambda}^0(M,\fg)\oplus\Omega_{\lambda}^2(M,\fg)$
and hence is Fredholm.  Since the dimension of $M$ is odd,
this operator has index zero.  If $A_0$ is a flat connection corresponding
to an isolated conjugacy class of 
irreducible representations of $\pi_1(M)$ in $G$,
the kernel of $\delta\cF_{s,A_0}|_{(A_0,0)}$ is zero; it follows that the cokernel
is also zero, so that $\delta\cF_{s,A_0}$ is an isomorphism at $(A_0,0).$
\end{proof}

\comment{check the transversality, also whether the statment of 
conj 3.19 and material above it should be modified in view of
the condition on A_0.}

As for the Wess-Zumino-Landau-Ginzburg model we have the following conjecture,
which may be viewed as giving a stochastic version of the Uhlenbeck compactness
theorem for solutions of the Yang-Mills equations.

\begin{Conjecture}\labell{conv3dgt}
The measure $\mu_{M,G}^{s,t,A_0}$ satisfies $$ \mu_{M,G}^{s,t,A_0}(1) < \infty$$
for $t$ sufficiently large.
\end{Conjecture}

Suppose $A_0$ is
a flat connection that corresponds
to an isolated conjugacy class of
irreducible representations of $\pi_1(M)$ in $G.$ 
Conjecture \ref{conv3dgt}
implies that for $t$ sufficiently
large, the measure $\mu_{M,G}^{s,t,A_0}$ is (up to normalization)
a Borel probability
measure on $\cA_{\lambda+2s+1} \oplus \Omega_{\lambda+2s+1}^3(M,\fg).$

\subsubsection{The phase of the determinant of $\delta \cF_{s,A_0}$}
The derivative $\delta \cF_{s,A_0}|_{(A,\xi)}$ is given by the Fredholm operator
$$D_{(A,\xi)}(\alpha,\eta)=(-\Delta_{A_0}+1)^s
(*d_{A_0}*\alpha+ d_A\alpha + *d_A*\eta + *[\alpha,*\xi]).$$
Suppose that $(A',0)$ is a regular point of $\cF_{s,A_0}.$
Then $D_{(A',0)}$ is invertible, and 
$(D_{(A',0)})^{-1} D_{(A,\xi)} = 1 + K((A',0);(A,\xi)),$
where
 $$K((A',0);(A,\xi)):
\Omega^1_{\lambda+2s+1}(M,\fg) \oplus \Omega_{\lambda+2s+1}^3(M,\fg)
\to
\Omega^1_{\lambda+2s+1}(M,\fg) \oplus \Omega_{\lambda+2s+1}^3(M,\fg)$$
is given by 
$$K((A',0);(A,\xi))(\alpha,\eta) =
(D_{(A',0)})^{-1}(-\Delta_{A_0}+1)^s
([A-A',\alpha] + *[A-A',*\eta] + *[\alpha,*\xi])
.$$

A similar computation to that done for the Wess-Zumino-Landau-Ginzburg and
Wess-Zumino-Novikov-Witten models shows that
$K((A',0);(A,\xi)) \in \cS_4(\Omega^1_{\lambda+2s+1}(M,\fg) \oplus \Omega_{\lambda+2s+1}^3(M,\fg))$
for all $(A,\xi).$
Thus the regularized determinant $det_4(1+K((A',0);(A,\xi)))$
is well defined,
and is nonzero at all regular points $(A,\xi)$
of $\cF_{s,A_0}.$  Let $U$ be the set of regular points of $\cF_{s,A_0}.$
Define $\Psi((A',0);\cdot): U \to S^1$ by
\begin{equation}\labell{det3dgt}
\Psi ((A',0);(A,\xi)) =
\frac{det_4(1+K((A',0);(A,\xi)))}{|det_4(1+K((A',0);(A,\xi)))|}.
\end{equation}

Then the analog of Conjecture \ref{wzlgconj} is the following.

\begin{Conjecture}\labell{casson}
Assume that Conjecture \ref{conv3dgt} holds.
Define the partition function $Z_s^t(M,G)$ by 
$Z_s^t(M,G) := \mu_{M,G}^{s,t,A_0}(\Psi((A',0),\cdot)).$  Then for
$t$ sufficiently large, $Z_s^t(M,G)$ is
independent of $t,$
and $|Z_s^t(M,G)|$ is a topological invariant of $M.$
If $M$ is a homology three-sphere and $A'=A_0$ is chosen
to be a flat connection corresponding to a conjugacy class
of isolated irreducible representations of $\pi_1(M)$ in $SU(2),$ then
$|Z_s^t(M,SU(2))|= |\chi(M)|$ where $\chi(M)$ is the Casson invariant of $M.$
\end{Conjecture}

Morally this conjecture says that for $t$ large, the measure
$\mu_{M,G}^{s,t,A_0}$ should be concentrated at the
transversal solutions of the nonlinear equation $\cF_{s,A_0} = 0,$
(or, equivalently, $\cF_{0,A_0} = 0$), which correspond to
conjugacy classes of isolated irreducible
representations of $\pi_1(M)$ in $G.$ These are the
representations which enter into the definition of 
the Casson invariant (see \cite{am}).

\begin{Remark} 
As in the case of the Wess-Zumino-Novikov-Witten model,
we expect that a nontrivial
phase of the determinant arises from the correction terms
in the regularized determinant.  Since the regularized determinant
appearing in (\ref{det3dgt}) is $det_4,$ equation (\ref{trcl}) leads
us to expect this "anomaly"
in $\Psi(A)$
to be the exponential of a cubic polynomial in $A.$  It seems reasonable to conjecture
that this cubic polynomial
is related to the Chern-Simons invariant of $A.$
\end{Remark}

\begin{Remark}\labell{pc}
Note that if one attempts to write down the
formal analog of the integral of the density (\ref{pb}), one obtains
a formal expression of the type 
\begin{equation}\labell{formal}
\int_\cA \exp(-S(A)) {\rm det} (\delta \cF_{s,A_0}),
\end{equation}

\noindent where $S(A) = |(-\Delta_{A_0} + 1)^s F_A|^2 + \dots.$

If one attempts to interpret the integral (\ref{formal}) as a perturbation of
a Gaussian measure, writing $S(A) = |(-\Delta_{A_0} + 1) d A|^2
+ I(A),$ a direct computation of Feynman diagrams
shows that the expectation of $|I(A)|^2$ in the appropriate
Gaussian measure diverges.
In the language of the physics literature, the integral (\ref{formal})
is {\em not} regularized.\footnote{Naive power counting indicates that
the expectation of $|I(A)|^2$ converges.  This does not imply convergence since
naive power counting is ineffective for theories with derivative
interactions.  One
example is the
expectation in the Gaussian measure with covariance
$(-\Delta_{A_0} + 1)^{-(2s + 1)}$ of the square of the cubic
vertex $ |<((-\Delta_{A_0} + 1)^s dA),((-\Delta + 1)^s (A^2))>|^2,$
which is divergent but for which naive power counting
predicts convergence.} Our method of constructing measures on these spaces
therefore differs in an essential way from standard renormalization
theory.\footnote{Readers familiar with the physics literature
might be surprised by the appearance of the non-gauge-invariant
Lagrangian (\ref{formal}) and by the absence of Faddeev-Popov
ghost terms.  The first problem may be remedied by replacing
the map $\cF_{s,A_0}$ by the map $\tilde{\cF}_{s,A_0}$ given by
\begin{equation}
\tilde{\cF}_{s,A_0}(A,\xi) = (-\Delta_{A} + 1)^s (*d_{A_0}*(A-A_0) + F_A + *d_A*\xi).
\end{equation}
which has the same zeros.  However, since this map also
corresponds to a gauge-fixed theory, there is little
to choose between the two maps, and we have preferred the simpler
map $\cF_{s,A_0}.$  The absence of Faddeev-Popov ghosts is
due to the cancellation of the ghost determinants between
bosons and fermions in a supersymmetric theory.
}
\end{Remark}

\subsection{Gauge theory in two and four dimensions}
\subsubsection{Two dimensional gauge theory}
Let $\Sigma$ be a 
Riemannian two-manifold, let $G$ be a compact Lie group.  Choose
an invariant Hermitian inner product on $\fg.$ 
Let
$\lambda < -1, $ let $s > -\lambda/ 2, $ and let 
$\cA_{\lambda+2s+ 1}$ be the space of connections on the trivialized principal
$G$-bundle over $\Sigma$  of Sobolev class $\lambda + 2s+ 1.$
Let $A_0$ denote the product
connection.  Let
$\Omega^i_\nu(\Sigma,\fg)$ denote the space of de Rham forms
on $\Sigma$ of Sobolev class $\nu$  with values in $\fg,$
and let $\star$
denote the Hodge star operator on differential forms
with values in $\fg$ obtained from the Riemannian
pairing on differential forms
and the hermitian inner product on $\fg.$
Let $E$ be a finite-dimensional vector space of dimension
$6g-6$, and let $f: \cA_{\lambda+ 2s+ 1} \to E$ be a map. Let 
$\cF_s : \cA_{\lambda+ 2s+ 1}  \to
\Omega^0_\lambda(\Sigma,\fg) \oplus \Omega^2_\lambda(\Sigma,\fg) \oplus E$
be given by 

\begin{equation}
\cF_s(A) = (-\Delta_{A_0} + 1)^s (F_A \oplus \star d_{A}\star (A-A_0)) \oplus f.
\end{equation}

Again, using Proposition \ref{sob}, $\cF_s$ is
a differentiable map of Banach manifolds, and its derivative is a Fredholm
map.  The index of this map depends on $f$; for $f =0 $ this index
is $6g-6,$ so we may expect the index of $\cF_s$ for generic $f$ to
be zero, and so to yield a measure on the space of connections
$\cA_{\l + 2s + 1}$.  The zeros of $\cF_s$ correspond to flat connections on $\Sigma$
satisfying conditions given by the map $f$ and so, for appropriate
choices of $E$ and $f,$
the resulting measure may be expected to
be related to intersection numbers on a moduli space of vector bundles.

\subsubsection{Four dimensional gauge theory}
Let $M$ be a 
Riemannian four-manifold, let $G$ be a compact Lie group.
Let $P$ be a principal $G$-bundle over $M.$  Choose an invariant Hermitian
metric on ${\rm ad}(P).$
Let $\lambda < -2, $ let $s > (-\lambda+1)/ 2, $ and let 
$\cA_{\lambda+2s+ 1}$ be the space of connections on
$P$ of Sobolev class $\lambda + 2s+ 1.$
Let $A_0$ be a smooth point of $\cA_{\lambda+2s+ 1}.$
Let
$\Omega^i_\nu(M,{\rm ad}(P))$ denote the space of de Rham forms
on $M$ of Sobolev class $\nu$  and values in ${\rm ad}(P).$
Let $\star$ denote the
Hodge star operator on ${\rm ad}(P)$-valued differential forms and let 
$\Omega^2_{\nu,+}(M,{\rm ad}(P))$ denote the space of self-dual
de Rham forms of Sobolev class $\nu.$
Let $E$ be a finite-dimensional vector space of dimension
$8p_1(P) - 3(b_2^+(M)-b_1(M) + 1)$, and let $f: \cA_{\lambda+ 2s+ 1} \to E$ be a map. Let 
$\cF_s : \cA_{\lambda+ 2s+ 1}  \to
\Omega^0_\lambda(M,{\rm ad}(P)) \oplus \Omega^2_{\lambda,+}(M,{\rm ad}(P)) \oplus E$
be given by 

\begin{equation}\labell{four}
\cF_s(A) = (-\Delta_{A_0} + 1)^s (F_A^+ \oplus \star d_{A}\star (A-A_0)) \oplus f,
\end{equation}

\noindent where $F_A^+$ is the self-dual part of the curvature of $A.$

Again, using Proposition \ref{sob} one can show that $\cF_s$ is
a differentiable map of Banach manifolds.  Its derivative is a Fredholm
map, the composition of $(-\Delta_{A_0} + 1)^s$ with a compact
perturbation of the standard elliptic operator $d_A^+ + \star d_A \star:
\Omega^1_{\l + 2s + 1}(M,{\rm ad}(P)) \to
\Omega^0_{\lambda+2s}(M,{\rm ad}(P)) \oplus
\Omega^2_{\lambda+2s,+}(M,{\rm ad}(P)).$

Again, for generic $f,$ we expect the index of $\delta\cF_s$
to be zero, and hence
to yield a measure on the space of connections
$\cA_{\l+ 2s +1}$.  The zeros of $\cF_s$ correspond to self-dual connections
on $M$
satisfying conditions given by the map $f$ and so, for appropriate
choices of $E$ and $f,$
the resulting measure may be expected to
be related to Donaldson invariants of $M.$

\begin{Remark}
The formal analogy between the path integrals appearing in the physics
literature on gauge
theory and the Mathai-Quillen formula was studied in \cite{aj}.
\end{Remark}
\begin{Remark}
There are other possible choices of the map $\cF_s$
in this example.  One alternative to the definition (\ref{four}) is
\begin{equation}
\tilde{\cF_s}(A) = (-\Delta_{A} + 1)^s (F_A^+ \oplus \star d_{A}\star (A-A_0)) \oplus f.
\end{equation}

This choice makes contact with the work of \cite{aj} as follows.
Write $\cA$ for the space of connections on $P$ and $\cG = Aut(P).$
In the spirit of \cite{aj}, the formal path integrals of \cite{tqft} morally
compute the Euler number of the vector bundle\footnote{We follow \cite{aj}
and ignore the locus where the action of $\cG$ is not free.}
${\mathcal V}:=\Omega^2_+(ad(P)) \times_\cG \cA \to \cA/\cG$
by a formula involving the 
section $A \to F_A^+,$ corresponding to the map $\tilde{\cF_0}.$  The same computation
applied to the section
$A \to (-\Delta_{A} + 1)^s (F_A^+),$
corresponding to the map $\tilde{\cF_s},$ should morally give the same result,
as it is a section of the same vector bundle.  In algebraic geometry there
is the familiar idea
of computing the number of zeros of a section of a vector bundle
by finding a different section of the same vector bundle where the computation
is simpler.  In our setting a different choice of section takes us
from the realm of formal path integrals to that of well-defined measures!

It is possible to give a construction similar to the one we have given
in this paper for measures corresponding to sections of the vector bundle
${\mathcal V} \to \cA/\cG$ rather than for
maps of Banach spaces.  For this purpose the map
$\tilde{\cF_s}$ (rather than $\cF_s$) is crucial, since it gives rise
to a section of the appropriate vector bundle.
Similar constructions exist for gauge theory in dimensions $2$ and $3,$
and indeed in higher dimensions.
Details will appear elsewhere.
\end{Remark}

\subsection{Quantum cohomology}  Let $\Sigma$ be a K\"ahler two-manifold,
and let $M$ be a K\"ahler manifold; denote by $T\Sigma$ and $TM$
the holomorphic
tangent bundles of $\Sigma$ and $M,$ respectively.
Let $\lambda < -1$ and $s > (-\lambda+1)/2.$  Let
${\rm Map}_{\lambda + 2s+ 1}(\Sigma, M)$ be the space of maps
from $\Sigma$ to $M$ of Sobolev class $\l + 2s +1.$
Let $V$ be the vector bundle
over ${\rm Map}_{\lambda + 2s+ 1}(\Sigma, M)$ whose fibre
at $\phi\in {\rm Map}_{\lambda + 2s+ 1}(\Sigma, M)$ is the space
$\Gamma_\lambda(\phi^*T^*M \otimes {T}\Sigma)$
of sections of $\phi^*T^*M\otimes {T}\Sigma$ of
Sobolev class $\lambda.$  Let
$\cF_s$ be the section of $V$ be given by 
$$\cF_s (\phi) = (-\Delta + 1)^s\bar{\partial} \phi.$$
In a sufficiently small neighborhood of any map $\phi$,  the
vector bundle $V$ can be trivialized as a space of maps into
a vector spaces.  The section $\cF_s$ can then be used to push forward
the white noise measure on these vector spaces to a measure on 
${\rm Map}_{\lambda + 2s+ 1}(\Sigma, M).$  The analog of the Mathai-Quillen
formula (\ref{mqf}) then differs from this measure by the usual phase
of the determinant as well as by an infinite series analogous
to the series $C_k$ appearing in (\ref{mqf}); this infinite
series can be shown to converge.  Details will appear elsewhere.
\begin{Remark}
The analogy between path integrals associated to quantum cohomology
in the physics literature and the Mathai-Quillen
formula associated formally to $\cF_0$ was studied in \cite{wu}.
\end{Remark}

\section{The Wess-Zumino-Novikov-Witten models on $\P^1.$}\labell{wznwsec}

In this section we make a more detailed study of the measure
$\mu_\cG^{s,t}$ constructed in Theorem \ref{wznwthm}.  First we prove
that $\mu_\cG^{s,t}$ is a Borel probability measure on $\cG.$  Since
$\mu_\cG^{s,t}(\cG)=1,$ there is no analog for $\mu_\cG^{s,t}$
of Conjecture \ref{wzlgconj} or Conjecture \ref{casson}.  Instead
we compute more subtle quantities given by expectation values of
functions in the measure $\mu_\cG^{s,t}.$  We show that these
expectation values give formulas which agree with formulas discovered
by Frenkel and Zhu \cite{fz} in the context of the theory of
vertex operator algebras.

\subsection{The measure $\mu_\cG^{s,t}$ is a probability measure.}
Recall that we are working with the a subgroup $G^\C \subset SL(n,\C),$ and
this identification
induces a hermitian inner product on $\fg^\C = {\rm Lie}(G^\C).$
Choose $\lambda < -2$ and $s > (-\lambda+1)/2.$
Recall that $\cG : = {\rm Map}_{*,\lambda + 2s + 1}(\P^1,G^\C)$
is the space of based maps from $\P^1$ to $G^\C$ of
Sobolev class $\lambda + 2s +1.$  Let 
$\tcg :={\rm Map}_{\lambda + 2s + 1}(\P^1,G^\C)$ denote the space
of all maps from $\P^1$ to $G^\C$ of Sobolev class $\l + 2s +1;$ like $\cG,$
this is a Lie group with the group structure given by pointwise multiplication.
Let $\cB := \Omega^{0,1}_{\lambda}(\P^1,\fg^\C)$ be the space of anti-holomorphic
one forms on $\P^1$ of Sobolev class $\lambda,$  and
let $\cA_{\lambda+ 2s}:= \Omega^{0,1}_{\lambda+ 2s}(\P^1,\fg^\C).$
The group $\tcg$ acts on $\cA_{\lambda + 2s}$ by the gauge action
$$g \cdot A = g A g^{-1} + \bar{\partial}g g^{-1}$$ for $g \in \tcg, A \in \cA_{\l+ 2s}.$
Recall that
$\cB = \G_\l(W)$ where $W = T^*\P^1 \otimes \fg^\C,$ and that
$\cB$ is therefore equipped with a white noise measure
$\mu_W.$

Recall that the map $\cF_s : \cG \to \cB$ is given by
$\cF_s (g) = (-\Delta + 1)^s  \bar{\partial} g g^{-1}.$
The map
$\cF_s$ is the composite of the map $\cF: \cG \to \cA_{\l+ 2s}$
given by $\cF(g) = \bar{\partial}g g^{-1}$ with the
isomorphism $(-\Delta + 1)^s:\cA_{\l + 2s} \to \cB.$
Since $\cF:\cG\to \cA_{\l+ 2s}$
gives the orbit of $0 \in \cA_{\l+ 2s} \in \cA,$
and since $\cG$ acts freely on $\cA_{\l+ 2s},$
$\cF_s$ is injective.  To determine
the complement of the image of $\cF_s,$ we make use of the 
stratification of the space $\cA_{\l + 2s}$ given
in Atiyah-Bott \cite{ab}, referring to the
work of Shatz \cite{shatz}, Narasimhan-Seshadri \cite{ns},
and Harder-Narasimhan \cite{hn}. See also 
\cite{groth1,atiyah,don1,dask,rade}.

\begin{Proposition}\labell{nsthm}
The image of $\cF_s$ consists of the complement of a countable union
of smooth, locally closed subvarieties $V_i$ of $\cB.$
Each of these subvarieties is of the form
$V_i = (-\Delta + 1)^s \tcg \cdot A_i$
where $A_i$ is a smooth element of $\Omega^{0,1}(\P^1,\fg^\C)$.
The subvariety $V_i$ has positive codimension in $\cB,$
and the normals to $V_i$ at $(-\Delta + 1)^s A_i$ are smooth.
\end{Proposition}

\begin{proof}
Since the map $(-\Delta + 1)^s : \cA_{\l + 2s} \to \cB$ is an 
isomorphism, it is enough to prove an analogous result for
the image of the map $\cF:\cG \to \cA_{\l + 2s}.$  This the orbit
of the point $0 \in \cA_{\l + 2s}$ under the action of $\tcg.$

The results of \cite{ab} show that the space $\cA_{\l+ 2s}$
can be stratified as $\cA_{\l+2s}= \coprod_\mu \cA_{(\mu)},$
where the $\cA_{(\mu)}$ are a countable family of $\tcg$ orbits
which are smooth locally closed subvarieties of $\cA_{\l+2s},$
each of finite codimension in $\cA_{\l+2s}.$  
On $\P^1$ these strata
are given by $\tcg$-orbits of a countable collection of 
elements $A_\mu \in \cA_{\l +2s}.$  (See
\cite{ab}, Section 14, pp. 608-609, \cite{dask}, or \cite{rade}, Theorem 2 and
its Corollary).

The $\tcg$-orbit of the trivial connection $A_0$ corresponding
under our identification to $0 \in \Omega^{0,1}_{\l+ 2s}(\P^1,\fg^\C)$
is the image of the map $\cF.$  We denote the corresponding
stratum by $\cA_{(0)}.$

Each stratum $\cA_{(\mu)}$ contains a smooth point by Lemma 14.8
of \cite{ab}.  The normals to $\cA_{(\mu)}$ at a smooth point
$A \in \cA_{(\mu)}$ are given by the kernel of the adjoint $D_A^*$
of the elliptic operator
$D_A:=\bar{\partial}+A:\Omega^{0}(\P^1,\fg^\C) \to \Omega^{0,1}(\P^1,\fg^\C),$
and are therefore smooth (see \cite{ab}, Section 14, p. 608).
The kernel of $D_A^*$ at smooth points was computed in
\cite{ab}, where it is shown to be zero for the stratum
$\cA_{(0)}$ and nonzero for every other stratum;
see Proposition 5.4, (5.10), (7.15), and (10.7) of \cite{ab}.

In the Yang-Mills context each of these strata is the
orbit under the complexified gauge group
of a solution of the Yang-Mills equation.
The minimum of the Yang-Mills functional, given by the trivial
connection, has the open orbit $\cF(\cG).$
\end{proof}

We therefore obtain the following result.

\begin{Proposition} \labell{embprop}
The map $\cF_s$ is a diffeomorphism of $\cG$ with the complement of 
a set of measure zero with respect to the white noise measure
$\mu_W^t$ on $\cB$ for any $t.$
\end{Proposition}
\begin{proof} 
Since the map $(-\Delta + 1)^s:\cA_{\l + 2s} \to\cB$ is an isomorphism, it suffices
to show that the subvariety $Z_i\subset \cA_{\lambda + 2s}$ given by 
$Z_i = \tcg \cdot A_i$ has measure
zero with respect to the Gaussian measure $\mu_{W,s}^t$
on the space $\cA_{\l+2s}=\Omega^{0,1}_{\lambda + 2s}(\P^1,\fg^\C).$

Let $a \in Z_i.$  By the inverse function theorem, we can
find $\epsilon > 0$ such that 
the orthogonal projection $B_\epsilon(a) \cap Z_i \to {T Z_i|_q} $
is an isomorphism for all $q \in B_\epsilon(a).$\footnote
{Given a subspace $T \subset \cA_{\l + 2s},$ denote by
$\pi_T: \cA_{\l +2s} \to T$ the orthogonal projection onto
$T.$  Define $f: Z_i \times {Z_i} \to T{Z_i}|_a \times T{Z_i}|_a$ by
$$f(p,q) = (\pi_{T{Z_i}|_a} \circ \pi_{T{Z_i}|_q} p, \pi_{T{Z_i}|_a} q).$$
Then $\delta f_{(a,a)}$ is an isomorphism.  By the inverse
function theorem, $f$ is a diffeomorphism in a neighborhood
of $(a,a).$  It follows that the orthogonal projection 
$\pi_{T{Z_i}|_q}|_{Z_i}$ is an isomorphism near $a$ for $q$ near $a.$}

Let $U=B_{\epsilon/2}(a).$  Since $Z_i = \tcg \cdot A_i,$ the point $a$ is
given by $a = g\cdot A_i$
for some $g \in \tcg.$  Let $g' \in \tcg$ be a smooth element of
$\tcg$ sufficiently close to $g$ so that $g'\cdot A_i \in U \cap Z_i.$
Then $p:=g'\cdot A_i$ is smooth.  The choice of $\epsilon$ is such that
the orthogonal projection $B_\epsilon(a) \cap Z_i \to {T Z_i|_p}$
is an isomorphism.

Let $\nu \in TZ_i|_p^\perp.$  Then $\nu$ is smooth, and there
exists $\alpha > 0$ so that if
$\delta, \delta' \in (0,\alpha)$ and $\delta \neq \delta',$ then
$$(U \cap Z_i + \delta  \nu) \cap (U \cap Z_i + \delta' \nu) =\emptyset.$$

Suppose that $\mu_{W,s}^t(U \cap Z_i) = \eta >0.$
By the Cameron-Martin formula (Proposition \ref{CM}),
the Radon-Nikodym derivative 
$d (T_{c \nu})_*\mu_{W,s}^t/d\mu_{W,s}^t$ is bounded from below for $q \in U$
and $c > 0$ by 
$$\frac{d (T_{c \nu})_*\mu_{W,s}^t}{d\mu_{W,s}^t}(q)=
 \exp(ct^2 <(-\Delta + 1)^{2s} \nu, q >  -\frac12 c^2 t^2||\nu||^2_{2s}) 
\geq \exp (-ct^2 ||\nu||_{-\l + 2s} ||q||_{\l + 2s}-\frac12 c^2 t^2||\nu||^2_{2s}).$$

\noindent Since for $q \in U,$ $||q||_{\l + 2s} < ||a||_{\l + 2s} + \epsilon,$ 
there exists $\beta > 0$ such that for all $q \in U$ and
all $c \in (0,\beta),$
$$\frac{d (T_{c \nu})_*\mu_{W,s}^t}{d\mu_{W,s}^t} (q)> \frac12.$$
Then if $c \in (0,\beta),$
$$\mu_{W,s}^t(U\cap Z_i + c \nu) > \frac12 \eta.$$
Let $\gamma = {\rm min~}\{\alpha,\beta\}.$  Then
$$\mu_{W,s}^t(\coprod_{n=1}^\infty (Z_i \cap U + \frac{\gamma}{2^n} \nu))=\infty,$$
which is impossible.  So it must be that  $\mu_{W,s}^t(U\cap Z_i)=0.$

Thus every point $p \in Z_i$ has a neighborhood in $Z_i$ of measure zero.
Since the space $\cA_{\l + 2s}$ has a countable basis for its topology,
so does the subspace $Z_i.$
Thus $Z_i$ is Lindel\"of, and can be covered by countably many open sets
of measure zero.  It follows that $\mu_{W,s}^t(Z_i) = 0.$
\end{proof}

As a corollary we get the following theorem.

\begin{Theorem}\labell{emb}
The measure $\mu_\cG^{s,t}$ is a Borel probability measure
on $\cG$ for all $t.$
\end{Theorem}

\subsection{The Frenkel-Zhu formula}

We now focus on the case $G=SU(n),$ so that $G^\C = SL(n,\C).$
To compare expectations in the measure $\mu_\cG^{s,t}$ to results arising
from the theory of vertex operator algebras, we introduce some functions
on $\cG.$
Recall that $\Omega^{0,1}_{\lambda + 2s}(\P^1,\fg^\C) = 
\G_{\l+2s}(W)$ where $W = T^*\P^1\otimes \fg^\C.$  We have equipped
$\P^1$ with a Kahler metric,
and the Lie algebra $\fg^\C$ comes with
an inner product given by   $(a,b) \to {\rm tr~} a^*b.$
Let us denote by $\star : \Omega^{0,1}(\P^1) \to \Omega^{1,0}(\P^1)$
the Hodge star operator.  Also, let
$\cC_z : T^*\P^1|_z \otimes (T^*\P^1|_z)^* \to \C$ denote the contraction.

The
inner product on $\fg^\C$ along with the Riemannian
structure on $\P^1$ give rise to a pairing
$$\langle,\rangle: 
\Omega^{0,1}_{\lambda + 2s}(\P^1,\fg^\C) \times
\Omega^{0,1}_{\lambda + 2s}(\P^1,\fg^\C)\to
H_{\l+2s}(\P^1);$$
here we have used Proposition
\ref{sob} to see that the inner product of elements in 
$\Omega^{0,1}_{\lambda + 2s}(\P^1,\fg^\C)$ is in
$H_{\l+2s}(\P^1).$

Similarly if we are given any hermitian
matrix $x$ we obtain from the pairing on $\fg^\C$ given
by $(a,b) \to {\rm tr~} a^* x b$ a pairing

$$\langle,\rangle_x: 
\Omega^{0,1}_{\lambda + 2s}(\P^1,\fg^\C) \times
\Omega^{0,1}_{\lambda + 2s}(\P^1,\fg^\C)\to
H_{\l+2s}(\P^1).$$

Let $x\in \fg$ be a hermitian matrix and let $z\in \P^1.$
We define a function $\hat{F}_{x,z}: \cG \to \R$ by

\begin{equation}\labell{vo}
\hat{F}_{x,z}(g) = \langle \bar{\partial}g g^{-1},
\bar{\partial}g g^{-1}\rangle_x (z)
\end{equation}

We also define ${F}_{x,z}: \cG \to \R$ by 

$$F_{x,z} := \hat{F}_{x,z} - \mu_{\cG}^{s,1}(\hat{F}_{x,z}). $$

Similarly, let $v  \in \Omega^{0,1}(\P^1,\fg^\C),$ let $z \in \P^1,$ and
define $f_{v,z}: \cG \to \C$ by 
$$ f_{v,z}(g) =  \langle v, \bar{\partial}g g^{-1}\rangle (z).$$

We recall some definitions from \cite{fz}.  Given a finite set $A,$
let $\cP_1(A)$ denote the set of partitions of $A$ into cycles and
chains containing precisely one chain.  Now let $A=\{1,\dots,n\}.$
Given $\alpha \in \cP_1(A)$, let us
write $\alpha$ explicitly as a collection of $k$ cycles and one chain
of length $m:$
$$\alpha = \{(a(1,1),\dots,a(1,j_1)),\dots,(a(k,1),\dots,a(k,j_k)),
[b(1),\dots,b(m)]\}.$$

\noindent Let $C_s(z,z')$ denote the Green's kernel
of the self-adjoint operator $(-\Delta + 1)^{2s}$
on $\Omega^{0,1}_{0}(\P^1).$ 
Then $C_s(z,z') \in( T^*\P^1|_{z'} )^*\otimes T^*\P^1|_z.$

Let  $v,v'\in  \Omega^{0,1}(\P^1,\fg^\C).$
We may write $v = \phi \xi, v' = \phi' \xi'$ where
$\phi,\phi' \in C^\infty(\P^1,\fg^\C)$ and 
$\xi, \xi' \in \Omega^{0,1} (\P^1).$

Let $ x_1,\dots,x_n \in\fg$, and let 
$z_1,\dots,z_n,z,w \in \P^1.$  
We define
$$f_{v,v',\alpha;s}(z,w;z_1,\dots,z_n) = \prod_{p=1}^k
{\rm tr ~}(x_{a(p,1)}\dots x_{a(p,j_p)} ) \times {\rm tr~}( {\phi}(z)^* x_{b(1)}\dots x_{b(m)}\phi'(w)) \times$$
$$\prod_{p=1}^k
\cC_{z_{a(p,1)}} \dots \cC_{z_{a(p,j_p)}}
(C_s(z_{a(p,1)},z_{a(p,2)})\otimes\dots \otimes
C_s(z_{a(p,j_p)},z_{a(p,1)}))\times$$
$$\cC_z \cC_{z_{b(1)}} \dots \cC_{z_{b(m)}}\cC_w
((\star\xi)(z)\otimes C_{s}(z,z_{b(1)})\otimes C_s(z_{b(1)},z_{b(2)})\otimes
\dots C_s(z_{b(m-1)},z_{b(m)}) \otimes
C_{s} (z_{b(m)},w) \otimes\xi'(w))
$$

Then we have the following result.

\begin{Theorem}\labell{fzt}
Let $x_1,\dots,x_n \in \fg$ and let $z_1,\dots,z_n,z,w\in \P^1.$
Let  $v,v'\in  \Omega^{0,1}(\P^1,\fg^\C).$ 
Then

\begin{equation}\labell{fzf}
\mu_{\cG}^{s,1}({\bar f_{v,z}} f_{v',w} F_{x_1,z_1} \dots F_{x_n,z_n}) = 
\sum_{\alpha \in \cP_1(A)} f_{v,v',\alpha;s}(z,w;z_1,\dots,z_n).
\end{equation}

\end{Theorem}

\begin{proof} 
By Theorem \ref{embprop}, we may transfer the computation to the space
$\cB,$ equipped with its white noise measure $\mu_W,$ using the map $\cF_s.$
For any $x \in \fg, z \in \P^1,$ we have
$F_{x,z}= \cF_s^* H_{x,z}$
where $H_{x,z}: \cG \to \C$ is given by 
$H_{x,z} : =\hat{H}_{x,z}-\mu_W(\hat{H}_{x,z}),$
and where $\hat{H}_{x,z}:\cB \to \C$ is given by
\begin{equation}
\hat{H}_{x,z}(A) = \langle (-\Delta + 1)^{-s} A,
 (-\Delta +1)^{-s} A \rangle_x (z)
\end{equation}
\noindent for $A \in \cB.$

Similarly, if $v  \in \Omega^{0,1}(\P^1,\fg^\C)$ and $z \in \P^1,$
$f_{v,z} = \cF_s^*h_{v,z}$
where $h_{v,z}:\cB \to \C$ is given by 
$$ h_{v,z}(A) =  \langle v, (-\Delta +1)^{-s} A \rangle (z).$$
Thus
$$\mu_{\cG}^{s,1}({\bar f_{v,z}} f_{v',w} F_{x_1,z_1} \dots F_{x_n,z_n})
=
\mu_W({\bar h_{v,z}} h_{v',w} H_{x_1,z_1} \dots H_{x_n,z_n}).$$

The expectation $\mu_W({\bar h_{v,z}} h_{v',w} H_{x_1,z_1} \dots H_{x_n,z_n}),$
like that of any polynomial,
can be computed by using formula (\ref{char}).  Let
$\phi_i \in \Omega^{0,1}_{-\l}(\P^1,\fg^\C)$, $i=1,\dots,2k.$   Let $t_1,\dots,t_{2k}
\in \C$ and let $\zeta := \sum_i t_i \phi_i.$
For $i = 1,\dots,2k,$ let $L_{\phi_i}:\cB \to \C$ be given
by $L_{\phi_i}(A) = <\phi_i,A>.$
Applying formula (\ref{char}) to $E_\zeta$
gives a generating
function $\mu_W(E_\zeta)$ for the expectations of all polynomials
in the $L_{\phi_i}.$  Explicitly, we have Wick's theorem:
Let ${\mathcal Q}_{2k}$ be the set of partitions
of the finite set $\{1,\dots,2k\}$ into a collection
of pairs $(a,b)$ where $a \in \{1,\dots,k\}$ and $b\in\{k+1,\dots,2k\}.$  Then

\begin{equation}\labell{wick}
\mu_W( \prod_{i=1}^{k} L_{\phi_i} \prod_{i=k+1}^{2k}\bar{L}_{\phi_i}) = 
\sum_{q \in {\mathcal Q}_{2k}}
\prod_{(a,b)\in q}
<\phi_a,\phi_b>.
\end{equation}

Applying equation (\ref{wick}) to the polynomial
${\bar h_{v,z}} h_{v',w} H_{x_1,z_1} \dots H_{x_n,z_n}$
gives (\ref{fzf}).

\end{proof}

The formula on the right hand side of equation (\ref{fzf}), when
evaluated at $s = 0,$ is (up to a constant) the formula
given in Theorem 2.3.1 of Frenkel and Zhu \cite{fz}, who work with the theory
of vertex operator algebras.
In our context we have obtained a type of analytic
continuation of this formula from a {\em commutative} algebra of functions
on $\cG.$

\begin{Remark}\labell{phirmk} The Wess-Zumino-Novikov-Witten model as it appears in the physics literature
differs from the formal path integral given by using the formula (\ref{pb}) for the map
$\cF_0$
in two ways.  First, the maps appearing in the Wess-Zumino-Novikov-Witten
path integral in the physics literature are maps
with values in the compact Lie group $G$ (rather than $G^\C$);
second, the partition
function is given by the integral with respect to this measure of
a power of a certain function $\Phi$.  Regarding the first issue,
the work of Wendt \cite{wendt} shows
that the quantization of the loop group $LG$ formally corresponds
to a path integral over a space of maps into $G^\C$, so that
a path integral over maps into $G^\C$ may be the correct model
for the study of the quantization of $LG.$  Another way of seeing
this is to repeat our construction for manifolds with boundary,
using the results of Donaldson \cite{don} to substitute for
Proposition \ref{nsthm}.  Details will appear elsewhere.

We now discuss the role of the function $\Phi.$
Recall that the analogy with finite dimensions leads us to believe
that one should compute expectations in the measure
$\mu_{\cG}^{s,1}(\Psi \cdot),$ where $\Psi$ is the phase of the
regularized determinant calculated in (\ref{phasewznw}).
In fact the physics literature refers to a calculation
of an expectation of the type $\mu(\Psi \Phi^k\cdot),$ where $k \in \Z,$
$\mu$ is a morally a measure on the space
$\cG_u:={\rm Map}_*(\P^1,G),$ and the function $\Phi:\cG_u\to S^1$ is given by
$\Phi = \exp(\frac{i}{12\pi}\int_B tr (\gamma^{-1} d \gamma)^3),$ where $B$ is a
three-manifold with $\partial B = \P^1,$ and $\gamma$ is an arbitrary
extension of $g$ to $B.$

The function $\Phi$ can be interpreted as a flat section
of the restriction to the orbit $\cG_u\cdot 0 \subset \cA_{\l + 2s}$ 
of a hermitian line bundle with connection $\cL\to \cA_{\l + 2s}$
\cite{rsw} on the space $\cA_{\l + 2s}.$
This section cannot be extended to $\cG\cdot 0$ as a flat section since $\cL$ is not flat
on $\cG \cdot 0.$
However, since the space $\cA_{\l + 2s}$ is affine, one can define
a function $\hat{\Phi}:\cA_{\l + 2s} \to S^1$ by choosing a base point
$A_0 \in \cA_{\l + 2s}$ and letting $\hat{\Phi}(A)$ be the holonomy of the
connection on $\cL$ along the straight line path from $A_0$ to $A.$  Now the 
determinant line bundle \cite{q} on $\cA_{\l + 2s}$ is given by  $\cL^h,$
where $h$ is the dual Coxeter number of $G,$
and the curvature of this line bundle 
is given \cite{q} by a {\em quadratic} K\"ahler potential.
This means it is reasonable to expect that
the function $\hat{\Phi}$ is a the exponential
of the pull back to $\cG$ of a
quadratic function on $\cA_{\l + 2s}.$
In this picture the phase $\Psi$ of the determinant appearing in (\ref{phasewznw})
should be replaced by the
section $\hat{\Phi}^h$ of the determinant line bundle $\cL^h,$
which arises from a method of regularizing determinants different from that
of \cite{simon}.
A result similar to that of Theorem \ref{fzt} holds
also if the measure $\mu_{\cG}^{s,1}$
is replaced by the complex valued measure
$\mu_{\cG}^{s,1}(\hat{\Phi}^{k+h} \cdot)$
where $\hat{\Phi}(g) = \exp(iQ(\bar{\partial} g g^{-1}))$ and $Q$ is
any quadratic function on $\cA_{\l + 2s}.$

It is reasonable to expect that the function 
$\Psi: \cA_{\l + 2s} \to S^1$ 
defined in equation (\ref{phasewznw}) is also a section of $\cL^h.$
Note that if this is the
case, the ideas of Remark \ref{heuristic} would again lead us to believe
that $\Psi$ should be of the form
$\Psi(g) = \exp(iQ(\bar{\partial} g g^{-1}))$ where $Q$ is
a (possibly nonlocal) quadratic function on $\cA_{\l + 2s}.$
\end{Remark}

\section{Conclusions}
\subsection{Analytic continuation.}  The measures constructed on
Banach spaces in the Wess-Zumino-Novikov-Witten model
give analytic continuations of correlation
functions arising in the associated theory of vertex operator
algebras.
This is reminiscent of the case
of correlation functions in Minkowskian quantum field theory, which
arise as analytic continuations of Euclidean correlation functions
which can be computed by Euclidean path integrals.
\subsection{Higher derivatives} The path
integrals in the examples studied above
correspond to partial differential equations of higher
order than those arising in the problems studied by the physicists.
One may wonder if
Lagrangians with higher derivatives may be useful as models in physics.
\subsection{Supersymmetry}  Our construction uses supersymmetry in
an essential way.
From this point of view the supersymmetric models are
more fundamental than the bosonic models; indeed one could try to
construct the bosonic model by showing that the inverse of an
appropriate determinant is integrable with respect to the supersymmetric
measure.
\subsection{Quantum field measures are not perturbations of Gaussian
measures.}  The measures appearing in our work are not in general perturbations
of Gaussian measures; they are locally {\em pushforwards} of such
measures.  In particular, factors of the type $(-\Delta + 1)^s$, though
they resemble dimensional regularization, may not suffice to regularize
the theory.
In effect the singularities of quantum field
theory fall into two categories:  Singularities arising from attempting
to define nonlinear functions on spaces of distributions, which disappear
when higher-order partial differential equations are considered,
and singularities arising from the fact that a pullback
measure may not be absolutely continuous with respect to Gaussian measure
(as the Cameron-Martin theorem shows is the case for
translations), which should not be removed (see Remark \ref{pc}).
The physics literature makes no distinction between these two types
of singularities.  Therefore the predictions of standard renormalization
theory, that gauge theory in dimensions higher than four should not exist,
and that sigma-models should exist only for manifolds of positive
curvature, may not apply to supersymmetric models constructed by our
methods.  In fact gauge
theory in six dimensions \cite{rw,dt,t} and eight dimensions \cite{bks}
has been studied in the physics literature, and corresponds to
mathematical results, as is the case for quantum
cohomology for manifolds which do not have positive curvature.
\subsection{Anomalies.}  The physics literature
assigns a special role to the phase of the determinant of the tangential
operator $\delta\cF_s.$   The invariants of maps in the finite
dimensional case involve the phase of the determinant of the
tangential operator.  But in infinite dimensions, it is only
the regularized determinant that can be defined.  The construction
of this determinant involves arbitrary choices (in our case, a choice
of a base point) and its phase
may therefore not lead to a topological invariant, or else fail to
have the same symmetries as the map $\cF_s.$  In the physics
literature the appearance of such an ``anomaly'' is held to
preclude the existence of a well-defined theory.  The fact that
we obtain a canonical construction of a measure but that the
construction of the phase of the determinant involved arbitrary
choices is in line with physicists' intuition, which singled out
the importance of these anomaly terms.
\subsection{Noncommutative and commutative algebras.}  Comparison
with the Frenkel-Zhu formula shows that a result arising from a
noncommutative vertex operator algebra can be viewed in terms of
the path integral as coming from a commutative algebra.  The original
Feynman-Kac formula is another reflection of this principle.  It may
be hoped that a better understanding of the role of path integrals
may be helpful in interpreting algebraic structures
arising in mathematics motivated by quantum field theory.
For example, a construction of the Wess-Zumino-Novikov-Witten
model in genus one should give a proof of the Kac
character formula \cite{kac} as adumbrated in Wendt \cite{wendt}.
\subsection{Equivariance.}  In the presence of a group action it should
be possible to define equivariant versions of the measures ${\cF_s}^*\mu;$
this will involve a convergent infinite series of the type appearing
in (\ref{mqf}).   The equivariant integral should satisfy a
localization formula as in finite dimensions.  In the case of the
Wess-Zumino-Novikov-Witten model the equivariant extension was
written down by Wendt \cite{wendt}; it is closely related to
the coset models studied by Gawedzki and Kupiainen \cite{gk}.
\subsection{Interpretation of the partition function}
It would be interesting to have a geometric
interpretation of finite dimensional analogs of
the partition functions $\cF_s^*\mu(\Psi^k)$ analogous
to the interpretation of $\cF_s^*\mu(\Psi)$ as a degree in formula
(\ref{ls}).\\~\\

{\bf Acknowledgments:}  I am grateful to I. Frenkel, V. Guillemin,
A. Jaffe, V. Jones, S. Sternberg, D. Stroock, S. T. Yau, and
G. Zuckerman for many helpful comments,
discussions and suggestions.  I would particularly like
to thank I. Frenkel for bringing the work of Frenkel-Zhu \cite{fz}
and Wendt \cite{wendt} to my attention, D. Stroock for his generosity
in answering questions on the analytical aspects of this work, and
V. Guillemin for many productive suggestions and corrections, and
encouragement during the course of this work.

After completing this work I received a copy of the work of C. Taubes
\cite{taubes}, which explores a different way of constructing quantum
field measures on some of the spaces we have considered.  I hope that
a combination of our methods with those of \cite{taubes} will be useful
in reviving interest in mathematical approaches to quantum field theory.

\end{document}